\documentclass{article}

\usepackage[utf8]{inputenc}
\usepackage[margin=20truemm]{geometry}
\usepackage{type1cm}
\usepackage{empheq}
\usepackage{mathtools}
\usepackage{comment}
\usepackage{amsthm, amsfonts}
\usepackage[pdftex]{graphicx}
\usepackage{bxtexlogo}
\usepackage[all]{xy} 
\usepackage{tikz}
\usepackage{stmaryrd}
\usepackage{amsmath, amssymb}
\usepackage{afterpage}
\usepackage{mathrsfs}
\usepackage{ascmac}
\usepackage{bm}
\usepackage{url}
\usepackage{makeidx}
\usetikzlibrary{graphs,graphs.standard}
\usepackage{caption}
\usepackage{latexsym}
\usepackage{afterpage}
\usepackage{nccmath}
\usepackage{enumerate}

\theoremstyle{definition}
\newtheorem{thm}{Theorem}[section]
\newtheorem{lem}[thm]{Lemma}
\newtheorem{cor}[thm]{Corollary}

\newtheorem{conj}[thm]{Conjecture}
\newtheorem{define}[thm]{Definition}
\newtheorem{prob}[thm]{Problem}
\newtheorem{ques}[thm]{Question}

\theoremstyle{remark}

\theoremstyle{definition}
\newtheorem*{thm*}{Theorem}
\newtheorem*{lem*}{Lemma}
\newtheorem*{prop*}{Proposition}
\newtheorem*{cor*}{Corollary}
\newtheorem*{ex*}{Example}
\newtheorem*{caution*}{Caution}
\newtheorem*{prob*}{Problem}
\newtheorem*{ques*}{Question}

\theoremstyle{definition}
\newtheorem*{define*}{Example}

\theoremstyle{remark}
\newtheorem*{rem*}{Remark}

\theoremstyle{remark} 
\newtheorem*{prf}{Proof}

\newcommand{\cG}{\mathcal{G}}

\newcommand{\cM}{\mathcal{M}}

\newcommand{\cS}{\mathcal{S}}

\newcommand{\mN}{\mathbb{N}}
\newcommand{\mZ}{\mathbb{Z}}

\newcommand{\set}[1]{\left\{#1\right\}}

\newcommand{\p}{\prime}

\newcommand{\lt}{\left}
\newcommand{\rt}{\right}

\newcommand{\ra}{\rightarrow}

\newcommand{\gir}{\mathrm{gir}}
\renewcommand{\mod}{\mathrm{mod}}

\newcommand{\ez}{\mathrm{ez}}

\renewcommand{\qed}{\hfill\Box}

\begin{document}

\title{Upper Bounds of the Odd Chromatic Number of a Graph in terms of its Thickness}
\date{} 
\author{S.Kitano}

\maketitle

\section{Introduction}
All considered graphs in this paper are simple, finite and undirected. For a vertex $v$ in a graph $G$, the neighbor of $v$ is a vertex adjacent to $v$. The neighborhood $N(v)$ is the set of the neighbors of $v$. The cardinality of $N(v)$ is the degree of $v$ and is denoted by $d(v)$. The value $\delta(G)=\min\set{d(v):v\in V(G)}$ and $\Delta(G)=\max\set{d(v):v\in V(G)}$ are the minimum and maximum degree of $G$, respectively. The average degree $\tilde{d}_G:=\sum_{v\in V(G)}\frac{d(v)}{|V(G)|}=\frac{2|E(G)|}{|V(G)|}$ of $G$ is the average value of the degrees of the vertices in $G$. The length of a shortest cycle contained in a graph $G$ is the girth of $G$ and is denoted by $\gir(G)$. A proper vertex coloring $\varphi:V(G)\ra\mN$ of a graph $G$ is a map such that if $uv\in E(G)$, then $\varphi (u)\neq\varphi(v)$. A graph is $k$-colorable if it has a proper $k$-coloring, that is, a proper coloring $\varphi:V(G)\ra\mN$ with the property that $\lt|\varphi\lt(V(G)\rt)\rt|\leq k$. For a graph $G$, the minimum integer $k$ such that $G$ has a $k$-coloring is the chromatic number of $G$ and is denoted by $\chi(G)$.\\
\subsection{History and the Original Definitions of Odd colorings}
An odd coloring was introduced by Petruševski and Škrekovski \cite{parity}. Let $G$ be a graph and $\varphi$ be a proper coloring of $G$. We say a vertex $v$ satisfies the parity condition if there exists a color $c$ such that $\lt|\varphi^{-1}(c)\cap N(v)\rt|$ is odd. A proper coloring $\varphi$ is an odd coloring if each non-isolated vertex $v\in V(G)$ satisfies the parity condition. If proper $k$-coloring $\varphi$ is an odd coloring, then we say that $\varphi$ is an odd $k$-coloring. A graph $G$ is odd $k$-colorable if it has an odd $k$-coloring. The minimum integer $k$ such that $G$ has an odd $k$-coloring is the odd chromatic number of $G$ and is denoted by $\chi_o(G)$. For a vertex $v\in V(G)$, define the set of colors $L^\ast_\varphi(v)=\set{c\in\set{1,\cdots,k}:\lt|\varphi^{-1}(c)\cap N(v)\rt|\equiv1(\mod\ 2)}$. If $L^\ast_\varphi(v)\neq\emptyset$, then choose an arbitrary color in the set $L^\ast_\varphi(v)$ and denote this color by $\varphi^\ast(v)$. If $L^\ast_\varphi(v)=\emptyset$, then define $\varphi^\ast(v)=0$. For example, consider the graph $G$ shown in Fig.\ref{ex:odd coloring graph}. The odd chromatic number of this graph is $4$. Two colorings $\varphi_1$ and $\varphi_2$ shown in Table \ref{ex:odd coloring table} are proper colorings of $G$. The coloring $\varphi_1$ is an odd coloring but $\varphi_2$ is not. This is because the vertex $v_3$ does not satisfy the parity condition with respect to $\varphi_2$.\\
\begin{figure}[htbp]
\centering
\begin{tabular}{cc}
  \begin{minipage}{0.45\linewidth}
    \centering
    \begin{tikzpicture}[every node/.style={circle,draw}]
      \node (1) at (0.5877852522924731291687059546390727685976524376431459910722724807, -0.809016994374947424102293417182819058860154589902881431067724311) {$v_4$};
      \node (2) at (-0.951056516295153572116439333379382143405698634125750222447305644, 0.3090169943749474241022934171828190588601545899028814310677243113) {$v_6$};
      \node (3) at (0.9510565162951535721164393333793821434056986341257502224473056444,0.3090169943749474241022934171828190588601545899028814310677243113) {$v_3$};
      \node (4) at (-0.587785252292473129168705954639072768597652437643145991072272480,-0.809016994374947424102293417182819058860154589902881431067724311) {$v_5$};
      \node (5) at (0,1) {$v_2$};
      \node (6) at (0,2) {$v_1$};
      \node (7) at (-2,0) {$v_7$};
      \foreach \u \v in {1/3,3/5,5/2,2/4,4/1,5/6}
        \draw (\u) -- (\v);
    \end{tikzpicture}
    \caption{A graph $G$ of order $7$ and size $6$.}
    \label{ex:odd coloring graph}
  \end{minipage} &
  \begin{minipage}{0.45\linewidth}
    \centering
    \begin{tabular}{|c||c|c|c||c|c|c|}
\hline
        $v$  & $\varphi_1(v)$  &  $L_{\varphi_1}^\ast(v)$ & $\varphi_1^\ast(v)$ & $\varphi_2(v)$  &  $L_{\varphi_2}^\ast(v)$ & $\varphi_2^\ast(v)$ \\
        \hline \hline
        $v_1$ & 1 & $\set{2}$ & 2 & 1 & $\set{2}$ & 2 \\
        $v_2$ & 2 & $\set{1}$ & 1 & 2 & $\set{1}$ & 1 \\
        $v_3$ & 1 & $\set{2,3}$ & 3 & 1 & $\set{}$ & 0 \\
        $v_4$ & 3 & $\set{1,4}$ & 1 & 2 & $\set{1,3}$ & 1 \\
        $v_5$ & 4 & $\set{1,3}$ & 1 & 3 & $\set{1,2}$ & 2 \\
        $v_6$ & 1 & $\set{2,4}$ & 4 & 1 & $\set{2,3}$ & 2 \\
        $v_7$ & 1 & $\set{}$ & 0 & 1 & $\set{}$ & 0 \\
        \hline
    \end{tabular}
    \caption{Two proper colorings $\varphi_1$, $\varphi_2$ of $G$ in Fig. \ref{ex:odd coloring graph} and corresponding vertex labels $L_\varphi^\ast$ and $\varphi^\ast$.}
    \label{ex:odd coloring table}
  \end{minipage}
\end{tabular}
\end{figure}
\begin{comment}
This section we introduce some results of odd coloring. 
First, we introduce the odd chromatic numbers of some representative graphs are as follows:\\
\begin{thm}{(\cite{odd-rem})}
\begin{displaymath}
\chi_o(T)=\left\{
\begin{array}{l}
1\ (T=K_1)\\
2\ (T\ is\ odd\ graph)\\
3\ (otherwise)
\end{array}
\right.
\end{displaymath}
\begin{displaymath}
\chi_o(P_n)=\left\{
\begin{array}{l}
n\ (n\leq2)\\
3\ (otherwise)
\end{array}
\right.
\end{displaymath}
\begin{displaymath}
\chi_o(C_n)=\left\{
\begin{array}{l}
3\ (n\equiv0(mod.3))\\
4\ (n\not\equiv0(mod.3),n\neq5)\\
5\ (n=5)\\
\end{array}
\right.
\end{displaymath}
\end{thm}
\end{comment}
Odd colorings of planar graphs are activity studied \cite{odd-rem,odd-sp,odd-planar-8,parity,sparse}. We introduce a conjecture and related results in those papers.\\
\begin{conj}{(\cite{parity})}
Every planar graph is odd $5$-colorable.\\
\end{conj}
Because $\chi_o(C_5)=5$, the bound is strict if the statement is in fact true. Petruševski and Škrekovski proposed this conjecture and proved that every planar graph is odd $9$-colorable. Subsequently, this was improved by Petr and Portier.\\
\begin{thm}{(\cite{odd-planar-8})}\label{thm:odd planar 8}
Every planar graph is odd $8$-colorable.\\
\end{thm}
There exist also several studies on planar graphs with restricted girth \cite{odd-sp,sparse}. Cho et al. proved the following \cite{odd-sp}. \\
\begin{thm}{(\cite{odd-sp})}\label{thm:girth 5 6}
Every planar graph with girth at least $5$ is odd $6$-colorable.\\
\end{thm}
\subsection{Thickness of Graphs}
In the early 1960s, Selfridge discovered that the complement of a planar graph of order $11$ is not planar \cite{bisurvey}. In other words, $K_{11}$ can not be represented as the edge union of two planar graphs. Determining whether the complete graph $K_9$ can be represented as the edge union of two planar graphs was a challenging problem, but its impossibility was independently proven in \cite{nine-1} and \cite{nine-2}. Based on this problem, Tutte introduced the concept of the thickness of a graph, which is the minimum number of planar graphs to represent as the edge union \cite{bisurvey}.\\
\begin{define}
The thickness of a graph $G$, denoted by $\theta(G)$ is the minimum $t$ such that the following two properties hold:
\begin{enumerate}[(1)]
\item $E_1$, $\cdots$, $E_t$ are pairwise disjoint sets such that $E(G)=E_1\cup\cdots\cup E_t$.
\item The subgraphs $(V(G), E_1),\cdots,(V(G), E_t)$ of $G$ are all planar.
\end{enumerate}
\end{define}
If $\theta(G)\leq2$, then the graph $G$ is called a biplanar graph. The coloring problem of biplanar graphs was introduced by Ringel in 1959 and is referred to as the Earth-Moon problem with the analogy of coloring Earth's and the Moon's maps \cite{moon}. This is still an open problem and the current status is between $9$ and $12$.\\
\begin{prob}
Determine the value of $\max\set{\chi(G):\theta(G)=2}$.\\
\end{prob}
At the time when Ringel proposed this problem, it was known that $\chi(G)\leq12$ \cite{tgt-thickness}. Ringel conjectured that $\chi(G)\leq8$. However, Sulanke discovered that $C_5\lor K_6$ is a biplanar graph and $\chi(C_5\lor K_6)=9$. More generality, the following results have been obtained.\\
\begin{thm}{(\cite{tgt-thickness})}\label{thm:proper thickness}
Let $G$ be a graph.
\begin{enumerate}[(a)]
\item\ $\max\set{\chi(G):\theta(G)=1}=4$.
\item $9\leq\max\set{\chi(G):\theta(G)=2}\leq12$.
\item $6\theta(G)-2\leq\max\set{\chi(G):3\leq\theta(G)\leq7}\leq6\theta(G)$.
\item $6\theta(G)-3\leq\max\set{\chi(G):7<\theta(G)}\leq6\theta(G)$.
\end{enumerate}
\end{thm}
In this paper, we investigate the upper bounds of the odd chromatic number for each level of graph thickness. A graph $G$ has a restriction of the average degree $\tilde{d}_G$ regarding to its thickness. We prove Lemma \ref{lem:t-planar average degree} concerning the average degree as it is used in the proof of Theorem \ref{thm:minor all t-planar}.\\
\begin{lem}\label{lem:t-planar average degree}
Let $\tilde{d}_G$ be the average degree of a graph $G$, then $\tilde{d}_G<6\theta(G)$.\\
\end{lem}
\begin{prf}
Let $t=\theta(G)$. There exists an edge partition $E_1$, $\cdots$, $E_t$ such that every $G_i=(V(G),E_i)$ is planar for $1\leq i\leq t$. By the Euler's formula, the size and the order of each graph $G_i$ satisfy $|E(G_i)|\leq3|V(G_i)|-6=3|V(G)|-6$. Therefore, $|E(G)|=\sum_{i=1}^{t}|E(G_i)|\leq (3|V(G)|-6)t$. This can be reformulated as $\tilde{d}_G =\frac{2|E(G)|}{|V(G)|}\leq6t-\frac{12}{|V(G)|}$. $\qed$\\
\end{prf}
\section{Odd $k^+$-Critical Graphs}

\subsection{Two useful Lemmas}
It is well known that $\chi(H)\leq\chi(G)$ for every graph $G$ and its subgraph $H$. If an equality holds if and only if $H=G$, then the graph $G$ is said to be color-critical. Furthermore, if $G$ is color-critical and $\chi(G)=k$, then we say that $G$ is a $k$-critical graph. We extend this concept to the odd chromatic number of $G$. In general, $\chi_o(H)\leq\chi_o(G)$ does not hold. (For example, consider the graph $G$ shown in Fig.\ref{ex:odd coloring graph}. While the graph $G$ has the graph $C_5$ as a subgraph, $\chi_o(C_5)=5>4=\chi_o(G)$.) However, considering a similar notion for odd coloring can be useful in this paper. \\
\begin{define}
Let $k\geq2$ be an integer. We say that a graph $G$ is an odd $k^+$-critical graph if (1) $\chi_o(G)\geq k$ and (2) $\chi_o(H)<k$ for every proper subgraph $H$ of $G$. We denote the class of odd $k^+$-critical graphs by $\cS_k$.\\
\end{define}
For example, the graph $C_5$ is an odd $4^+$-critical graph because $\chi_o(C_5)=5$ and every proper subgraph of $C_5$ is a forest. Lemma \ref{lem:fujie} and Lemma \ref{lem:n-easy key lemma} refer to the structure of odd $c^+$-critical graphs. Lemma \ref{lem:fujie} asserts that there are no vertices of odd low degree and no edges connecting vertices of low degree.\\
\begin{lem}\label{lem:fujie}
Let $G$ be an odd $k^+$-critical graph, where $k\geq4$. Then the following conditions hold.
\begin{enumerate}[(1)]
\item For every vertex $v\in V(G)$, either $d(v)$ is even or $d(v)\geq\lt\lfloor\frac{k}{2}\rt\rfloor$.
\item For every edge $e=v_0v_1\in E(G)$, if $d(v_0)\equiv d(v_1)(\mod.2)$, then either $\max\set{d(v_0),d(v_1)}>\frac{k-1}{2}$ or $d(v_0)=d(v_1)=\frac{k-1}{2}$.
\end{enumerate}
\end{lem}
\begin{prf}
(1) Let $v$ be a vertex. If $d(v)$ is odd and $d(v)<\lfloor\frac{k}{2}\rfloor$, then the graph $G-v$ has an odd $(k-1)$-coloring $\varphi$. Define $C=\set{\varphi(w),\varphi^\ast(w)|w\in N(v)}\backslash\set{0}$. Since $|C|\leq2d(v)<k-1$, there exists a color $\alpha\not\in C$. Define $\varphi(v)=\alpha$ and obtain an odd $(k-1)$-coloring $\varphi$ of $G$. $\qed$\\
(2)Let $v_0v_1$ be an edge such that $d(v_0)\equiv d(v_1)(\mod\ 2)$. Assume that $d(v_0)\leq d(v_1)\leq\frac{k-1}{2}$ and $d(v_0)+d(v_1)<k$. We have $d(v_0)\equiv d(v_1)\equiv0(\mod.2)$ by (1). The graph $G-\set{v_0,v_1}$ has an odd $(k-1)$-coloring $\varphi$. For any $i=0,1$, there exists a color in $\varphi\lt(N(v_{1-i})\backslash\set{v_i}\rt)$ such that it appears an odd number of times in $N(v_{1-i})\backslash\set{v_i}$. We arbitrarily choose one and denote it by $\alpha_i$. Define $C_i=\lt(\set{\varphi(w),\varphi^\ast(w)|w\in N(v_i)\backslash\set{v_{1-i}}}\backslash\set{0}\rt)\cup\set{\alpha_{1-i}}$. By $|C_i|\leq2\lt(d(v_i)-1\rt)+1<k-1$, there exists a color $\beta_i\not\in C_i$. Define $\varphi(v_i)=\beta_i$ and obtain an odd $(k-1)$-coloring $\varphi$ of $G$. $\qed$\\
\end{prf}
In \cite{odd-sp}, a useful lemma was introduced in the proof of Theorem \ref{thm:girth 5 6}. This lemma is related to the concept of an easy vertex. Before stating the lemma, let us define some terms. A vertex of degree exactly (at most) $d$ is called a $d$-vertex ($d^-$-vertex). If a $d$-vertex $u$ is adjacent to a vertex $v$, we say that $u$ is a $d$-neighbor of $v$. Denote the set of $d$-neighbors of $v$ by $N_d(v)$. Similarly define a $d^-$-neighbor and $N_{d^-}(v)$. A vertex $v$ is said to be easy if the degree of $v$ is odd or $N(v)$ contains a $2$-vertex. For a vertex $v$, the set of easy vertices in $N(v)$ is denoted by $N_{\ez}(v)$.\\
\begin{lem}\label{lem:key lemma}{(\cite{odd-sp})}
Let $G$ be an odd $k^+$-critical graph, where $k\geq6$. The graph $G$ has no $1$-vertex and $2d(v)\geq|N_2(v)|+|N_{\ez}(v)|+k-1$ for every easy vertex $v$.\\
\end{lem}
We generalize this lemma. First, we extend the definition of an easy vertex.\\
\begin{define}
Let $n$ be a positive integer. A vertex $v$ is said to be $n$-easy if its degree is odd or $\bigcup_{i=1}^nN_{2i}(v)\neq\emptyset$. For a vertex $v$, the set of $n$-easy vertices in $N(v)$ is denoted by $N_{n-\ez}(v)$.\\
\end{define}
\begin{lem}\label{lem:n-easy key lemma}
Let $n$ be a positive integer. Let $G$ be an odd $k^+$-critical graph, where $k\geq4n+2$. Then $2d(v)\geq\lt(\sum_{i=1}^n|N_{2i}(v)|\rt)+|N_{n-\ez}(v)|+k-1$ for every $n$-easy vertex $v$.\\
\end{lem}
\begin{prf}
By Lemma \ref{lem:fujie}, the degree of each vertex is either even or at least $2n+1$. We prove the inequality by contradiction. For an $n$-easy vertex $v$, we assume that $2d(v)\leq|X|+|N_{n-\ez}(v)|+k-2$. Write $X=\bigcup_{i=1}^nN_{2i}(v)$ and $Y=N_{(2n+1)^+}(v)$. The graph $G-\lt(X\cup\set{v}\rt)$ has an odd $(k-1)$-coloring $\varphi_0:V(G-\lt(X\cup\set{v}\rt))\ra\set{1,\cdots,k-1}$. We extend $\varphi_0$ to a coloring $\varphi_1:V(G)\ra\set{1,\cdots,k}$. For any $x\in X$, define $\varphi_1(x)=k$. If $w\not\in X\cup\set{v}$, then define $\varphi_1(w)=\varphi_0(w)$. We first define the color $\varphi_1(v)$. Define $A=\varphi_0(Y)\cup\lt(\set{\varphi_0^\ast(y)|y\in Y\backslash N_{n-\ez}(v)}\backslash\set{0}\rt)$. There exists a color $\alpha\in \set{1,\cdots,k-1}\backslash A$, since $|A|\leq\lt(|Y|+|X|\rt)+|Y\backslash N_{n-\ez}(v)|\leq d(v)+\lt(d(v)-|X|-|N_{n-\ez}(v)|\rt)<k-1$. Define $\varphi_1(v)=\alpha$. From the above we obtain
\begin{displaymath}
\varphi_1(w)=\left\{
\begin{array}{ll}
\varphi_0(w)& (w\not\in X\cup\set{v})\\
k & (w\in X)\\
\alpha& (w=v).
\end{array}
\right.
\end{displaymath}
We next define a proper $(k-1)$-coloring $\varphi_2$ of $G$. If $w\not\in X$, define $\varphi_2(w)=\varphi_1(w)$. When $X$ is nonempty, say $X=\set{x_1,x_2,\cdots,x_l}$. Define $B_i=\lt(\set{\varphi_2(y),\varphi^\ast(y):y\in N(x_i)\backslash\set{v}}\backslash\set{0}\rt)\cup\set{\alpha}$. Because $|B_i|\leq2d(x_i)-1\leq4n-1<k-1$, there exists a color $\beta_i\in \set{1,\cdots,k-1}\backslash B_i$ for $1\leq i\leq l$. Define $\varphi_2(x_i)=\beta_i$ and obtain
\begin{displaymath}
\varphi_2(w)=\left\{
\begin{array}{ll}
\varphi_0(w)& (w\not\in X\cup\set{v})\\
\beta_i & (w=x_i)\\
\alpha& (w=v).
\end{array}
\right.
\end{displaymath}
Let $Z=\set{z\in V(G):L_{\varphi_2}^\ast(z)=\emptyset}$. If $Z=\emptyset$, then $\varphi_2$ is an odd $(k-1)$-coloring of $G$. Otherwise, write $Z=\set{z_1,z_2,\cdots,z_m}$ and $z_i$ is an $n$-easy vertex and even degree for $1\leq i\leq m$. For each $z_i$, select a vertex $u_i\in\bigcup_{j=1}^nN_{2j}(z_i)$. Let $\psi_0=\varphi_2$, $Z_0=Z$. We define $\psi_i$ by recoloring $\psi_{i-1}$. For $1\leq i\leq m$, define $C_i=\lt(\set{\psi_i(y),\psi_{i-1}^\ast(y):y\in N(u_i)}\backslash\set{0}\rt)\cup\set{\psi_i(u_i)}$. There exists a color $\gamma_i\not\in C_i$ since $|C_i|\leq2d(u_i)\leq4n<k-1$. We define $\psi_i$ as follows and $Z_i=\set{z\in V(G):L_{\psi_i}^\ast(z)=\emptyset}$ for $1\leq i\leq m$.
\begin{displaymath}
\psi_i(w)=\left\{
\begin{array}{ll}
\psi_{i-1}(w)& (w\in V(G)\backslash\set{u_i})\\
\gamma_i & (w=u_i).
\end{array}
\right.
\end{displaymath}
Since $Z_0\supset Z_1\supset\cdots\supset Z_m=\emptyset$, $\psi_m$ is an odd $(k-1)$-coloring of $G$. $\qed$\\
\end{prf}
\subsection{Odd $k^+$-Critical Graphs with Respect to minor}
In the previous section, we defined an odd $k^+$-critical graph with respect to taking subgraphs. In this section, we consider a similar concept by taking minors instead.\\
\begin{define}
Let $k\geq2$ be an integer. We say that a graph $G$ is an odd $k^+$-critical graph with respect to taking minors if (1) $\chi_o(G)\geq k$ and (2) $\chi_o(H)<k$ for every proper minor $H$ of $G$. Denote the class of odd $k^+$-critical graphs with respect to taking minors by $\cM_k$.\\
\end{define}
Since minor operations include taking subgraphs, $\cM_k\subset\cS_k$. Considering a more restricted graph class allows for stronger assertions. First, we introduce a lemma. Subsequently, we provide an upper bound of the odd chromatic number of a graph whose minors have the thickness at most $t$.\\
\begin{lem}\label{lem:minor delta}
Let $k\geq2$ be an integer. If $G\in\cM_k$, then $\delta(G)\geq\lt\lfloor\frac{k}{2}\rt\rfloor$.\\  
\end{lem}
\begin{prf}
Assume, to the contrary, that a graph $G\in\cM_k$ has a $\lt(\lfloor\frac{k}{2}\rfloor-1\rt)^-$-vertex $v$. The graph $G-v$ has an odd $(k-1)$-coloring $\varphi$. Define $C=\set{\varphi(w),\varphi^\ast(w):w\in N(v)}\backslash\set{0}$. Since $|C|=2d(v)<k-1$, there exists a color $\alpha\not\in C$. Define $\varphi(v)=\alpha$ and obtain an odd $(k-1)$-coloring $\varphi$.$\qed$\\
\end{prf}
\begin{thm}\label{thm:minor all t-planar}
For $t\geq1$, if $G$ is a graph whose minors have the thickness at most $t$, then $\chi_o(G)\leq12t-1$.\\
\end{thm}
\begin{prf}
Assume, to the contrary, that there exists a graph $G\in\cM_{12t}$ such that every minor has thickness at most $t$. By Lemmas \ref{lem:t-planar average degree} and \ref{lem:minor delta}, we have $6t\leq\delta(G)\leq\tilde{d}(G)<6t$. $\qed$\\
\end{prf}
\begin{comment}
The odd chromatic number of a complete k-partite graph $G=K_{a_1,\cdots,a_k}$ where the size of each part is even is $k+2$. Therefore, according to Theorem \ref{thm:minor all t-planar}, if $k\geq12t-2$, then $G$ has a thickness-$(t+1)$ minor.\\
\begin{cor}
The complete $(12t-2)$-partite graph $K_{2,\cdots,2}$  has a thickness-$(t+1)$ minor.\\
\end{cor}
\end{comment}
\section{The main result}
The graph $K_n^\ast$ is constructed by subdividing each edge of the complete graph $K_n$ of order $n$ exactly once, and it satisfies $\delta(K_n^\ast)=2$ and $\gir(K_n^\ast)=4$. For each $n$, the graph $K_n^\ast$ is a biplanar graph and $\chi_o(K_n^\ast)=n$. Hence, there is no upper bound of the odd chromatic numbers of the class of biplanar graphs. This is surprising. From Theorem \ref{thm:proper thickness}, the class of biplanar graphs has an upper bound of the chromatic number, and from Theorem \ref{thm:odd planar 8}, the class of planar graphs has an upper bound of the odd chromatic number. However, the class of biplanar graphs does not have an upper bound of the odd chromatic number. We prove that the odd chromatic numbers of biplanar graphs with some specific conditions can be bounded by a constant and generalize this result to graphs with thickness at most $t$. Our result requires that the girth and minimum degree are large enough. More specifically, for example, a biplanar graph $G$ with $\gir(G)\geq6$ and $\delta(G)\geq3$ is odd $12$-colorable.\\
\begin{thm}\label{thm:girth 6 delta 3 color 12}
If $G$ is a graph with $\delta(G)\geq2\theta(G)-1$ and $\gir(G)\geq6$, then $\chi_o(G)\leq6\theta(G)$.\\
\end{thm}
To provide a proof, we use the discharging method. This approach uses the sum of universal quantities determined from the Euler's formula. Denote the set of components of a planar graph $G$ by $C(G)$ and the set of face of the fixed plane embedding by $F(G)$. Define $X(G)=V(G)\cup F(G)\cup C(G)$ and the map $\mu_G:X(G)\ra\mZ$ by
\begin{displaymath}
\begin{array}{lll}
\mu_G(f)=d(f)-6&\lt(f\in F(G)\rt)\\
\mu_G(v)=2d(v)-6&\lt(v\in V(G)\rt)\\
\mu_G(c)=6&\lt(c\in C(G)\rt).
\end{array}{}
\end{displaymath}
Define $S=\sum_{x\in X(G)}\mu_G(x)$. We can easily compute $S=-6$ by Euler's formula. We now proceed to the proof of Theorem \ref{thm:girth 6 delta 3 color 12}.\\
\begin{prf}
Let $t\geq2$. Let $G\in\cS_{6t+1}$ be a graph with $\gir(G)\geq6$, $\delta(G)\geq2t-1$, and $\theta(G)\leq t$. There exists an edge partition into $t$ planar graphs denoted by $G_1$, $\cdots$, $G_t$. Let $d_i(v)$ be the degree of a vertex $v\in V(G)$ in $G_i$. Fix plane embeddings of $G_1$, $\cdots$, $G_t$. Let $\mu_i$ be the map $\mu_{G_i}:X(G_i)\ra\mZ$. Define $S=\sum_{i=1}^t\lt(\sum_{x\in X(G_i)}\mu_i(x)\rt)$ and obtain $S=-6t$. Let $\mu^\star_i(x)$ be the final charge of $x\in X(G)$ after the following discharging rule:\\

If $d(v)<3t$, then each vertex $w\in N(v)$ sends a charge of $1$ to the vertex $v$.\\

Define $S^\star=\sum_{i=1}^t\lt(\sum_{x\in X(G_i)}\mu^\star_i(x)\rt)$ and obtain $S^\star=-6t$ since $S^\star=S$. We prove the following two conditions.
\begin{enumerate}[(a)]
\item For each $x\in F(G_i)\cup C(G_i)$ satisfies $\mu^\star_i(x)\geq0$.
\item For each $v\in V(G)$ satisfies $\mu^\star(v):=\sum_{i=1}^t\mu^\star_i(v)\geq0$.
\end{enumerate}
For (a), let $f$ be a face in $G_i$. Since girth of $G$ is at least $6$, $\mu^\star_i(f)=\mu_i(f)\geq6-6=0$. For each $c\in C(G_i)$, clearly $\mu_i^\star(c)=\mu_i(c)=6$. We next verify (b). Let $v$ be a vertex in $V(G)$. Note that for any vertex $v$, it holds that $\mu(v):=\sum_{i=1}^t\mu_i(v)=2d(v)-6t$. If $d(v)$ is odd, then $d(v)\geq3t$ by Lemma \ref{lem:fujie}. Therefore, we only need to consider cases where $d(v)<6t$ and $\bigcup_{j=1}^{\lfloor\frac{3}{2}t\rfloor}N_{2j}(v)$ is non-empty. First, we consider the case where $2t\leq d(v)<3t$. If $\bigcup_{j=1}^{\lfloor\frac{3}{2}t\rfloor}N_{2j}(v)$ is non-empty, then $\lt|\bigcup_{j=1}^{\lfloor\frac{3}{2}t\rfloor}N_{2j}(v)\rt|\leq2d(v)-6t<0$ by Lemma \ref{lem:n-easy key lemma}. Hence, we obtain $\mu^\star(v)\geq3d(v)-6t\geq0$. Next, we consider the case where $3t\leq d(v)<6t$ and $\bigcup_{j=1}^{\lfloor\frac{3}{2}t\rfloor}N_{2j}(v)$ is non-empty. We have $\lt|\bigcup_{j=1}^{\lfloor\frac{3}{2}t\rfloor}N_{2j}(v)\rt|\leq2d(v)-6t$ by Lemma \ref{lem:n-easy key lemma}. So, we obtain $\mu^\star(v)\geq2d(v)-6t-\lt|\bigcup_{j=1}^{\lfloor\frac{3}{2}t\rfloor}N_{2j}(v)\rt|\geq0$.\\
Therefore, $S^\star\geq0$ is established, which contradicts our earlier calculation $S^\star=-6t$. Consequently, there exists no graph $G$ such that $\delta(G)\geq2t-1$, $\gir(G)\geq6$ and $\theta(G)\leq t$ in $\cS_{6t+1}$. Thus, Theorem \ref{thm:girth 6 delta 3 color 12} is obtained from the minimality of $G$. $\qed$\\
\end{prf}
If $G$ is $(2\theta(G)-1)$-connected, then $\delta(G)\geq2\theta(G)-1$.\\
\begin{cor}
If $G$ is a $(2\theta(G)-1)$-connected graph with $\gir(G)\geq6$, then $\chi_o(G)\leq6\theta(G)$.\\
\end{cor}
Regarding the proof of Theorem \ref{thm:girth 6 delta 3 color 12}, considering surplus charges allows us to establish the existence of several $(2\theta(G)-1)^-$-vertices.\\
\begin{cor}\label{cor:-vertices}
Let $n_1$ be any natural number, $n_2=3$ and $n_k=2$ for $k\geq3$. If $G$ is a graph with $\gir(G)\geq6$ such that the number of $(2\theta(G)-1)^-$-vertices is at most $n_{\theta(G)}$, then $\chi_o(G)\leq6\theta(G)$.\\
\end{cor} 
\begin{prf}
Let $\tau=\theta(G)\geq2$. If the graph $G$ considered in the proof of Theorem \ref{thm:girth 6 delta 3 color 12} contains a $(2\tau-1)^-$-vertex $y$, then $\mu^\star(y)=3d(y)-6\tau\geq-6(\tau-1)$. By (a), if $Y=\set{y\in V(G):2\leq d(y)\leq2\tau-1}$ is not empty, then $S^\star=\sum_{i=1}^\tau\lt(\sum_{x\in X(G_i)}\mu^\star_i(x)\rt)\geq-6|Y|(\tau-1)$. Hence, if $|Y|<\frac{2\tau}{\tau-1}$, then $\chi_o(G)\leq6\tau$ since $-6|Y|(\tau-1)>-6\tau$. The Corollary \ref{cor:-vertices} is correct since $n_\tau<\frac{2\tau}{\tau-1}$. $\qed$\\
\end{prf}
\section{Problems for Further Study}
We end this paper with a question and a conjecture. There may be potential for improving the upper bounds established by Theorems \ref{thm:minor all t-planar} and \ref{thm:girth 6 delta 3 color 12}.\\
\begin{ques}
Determine whether the upper bounds of Theorems \ref{thm:minor all t-planar} and \ref{thm:girth 6 delta 3 color 12} are strict.\\
\end{ques}
\begin{comment}
For every positive integer $t$, the set of thickness-$t$ graphs that its minors are all thickness-$t$ graph is minor closed, so there exists forbidden minor list.\\
\begin{ques}
Determinant the set of thickness-$t$ graphs that its minors are all thickness-$t$ graph.\\
\end{ques}
\end{comment}
We propose a conjecture. The class of biplanar graphs has no upper bound of the odd chromatic number because the graph $K_n^\ast$ is a biplanar graph and its odd chromatic number is $n$.\\
\begin{conj}
Let $n\geq9$ be an integer. If $G$ is a biplanar graph without $K_{n+1}^\ast$ as a subgraph, then $\chi_o(G)\leq n$.\\
\end{conj}

\end{document}